\documentclass{amsart}      
\usepackage{amsmath,amsfonts,amssymb}
\usepackage[dvips]{graphicx}
\usepackage{psfrag}
\usepackage[T1]{fontenc}
\usepackage[latin9]{inputenc}

\theoremstyle{definition} 
\newtheorem{defs}{Definition}[section]

\theoremstyle{plain} 

\newtheorem{lem}[defs]{Lemma}
\newtheorem{thm}[defs]{Theorem} 

\theoremstyle{remark} 
\newtheorem{remark}[defs]{Remark}
\newtheorem{Example}[defs]{Example}

\newcommand{\R}{\mathbb{R}}

\newcommand{\Q}{\mathbb{Q}}
\newcommand{\Z}{\mathbb{Z}}
\newcommand{\CC}{\mathcal{C}}
\newcommand{\HH}{\mathcal{H}}
\newcommand{\QQ}{\mathcal{Q}}

\begin{document}

\title[Ends of strata]{Ends of strata of the moduli space of quadratic differentials}

\author{Corentin Boissy
}


\address{ 
Aix Marseille Université, CNRS,  LATP, UMR 7353, 13453 Marseille France } 
             \email{corentin.boissy@latp.univ-mrs.fr}

\date{\today}

\maketitle

\begin{abstract}
Very few results are known about the topology of the strata of the moduli space of quadratic differentials. In this paper,  we prove  that any connected component of such strata has only one topological end.
 A typical flat surface in a neighborhood of the boundary is naturally split by a collection of parallel short saddle connections, but the discrete data associated to this  splitting can be quite difficult to describe. In order to bypass these difficulties, we use the Veech zippered rectangles construction and the corresponding (extended) Rauzy classes.

\end{abstract}

\section{Introduction}

We study compact surfaces endowed with a flat metric with isolated conical singularities and $\Z/ 2\Z$ linear holonomy. Such a surface is naturally identified with a Riemann surface endowed with a meromorphic quadratic differential with at most simple poles. The moduli space of such surfaces with fixed combinatorial data is a noncompact complex-analytic orbifold $\QQ$ and is called a stratum of the moduli space of quadratic differentials. 

There is an obvious way to leave any compact of $\QQ$ by rescaling the metric so that the area tends to infinity or to zero. Hence we usually consider normalized strata that correspond to area one flat surfaces. A normalized strata is still noncompact, and a neighborhood of the boundary corresponds to flat surfaces with a short geodesic joining two singularities (not necessary distinct).

Very few results are known about the topology of these strata. Kontsevich, Zorich and Lanneau have classified their connected components (see \cite{Kontsevich:Zorich} and \cite{Lanneau:cc}). Eskin, Masur and Zorich have described the neighborhood of the \emph{principal boundary}, which corresponds to the neighborhood of the boundary after removing a subset of negligible measure (see \cite{EMZ,Masur:Zorich}). For the special case of genus zero flat surfaces, we have proven that the corresponding strata have only one topological end (see \cite{B2}). In this paper, we extend this result to all strata.

\begin{thm}\label{end}
Let $\CC$ be a connected component of a stratum of the moduli space of quadratic differentials and let $\CC_1$ be the subset of $\CC$ that corresponds to area one flat surfaces.
Then, $\CC_1$ has only one topological end.
\end{thm}

The most natural approach is to describe a typical flat surface in the neighborhood of the boundary. A saddle connection is a geodesic joining two singularities. A flat surface is near the boundary if it has a saddle connection of short length. One can look at the set of saddle connections that are of minimal length. In general, there are several such saddle connections and we can show that they are parallel for a generic flat surface, and they stay parallel and of the same length for any small perturbation of the surface. One can associate to such a collection of saddle connections a ``configuration'' that describes how the collection splits the surface (see \cite{EMZ,Masur:Zorich}). 
The number of different configurations tends to infinity when the genus tends to infinity, and no canonical way is known to describe all the configurations associated to a connected component of a stratum (see \cite{Masur:Zorich,B1}). We show in the Appendix  some examples that illustrate the difficulties that appear using this direct approach.
\medskip

In order to bypass these difficulties, we use \emph{generalized permutations} that give another combinatorial description of a surface in a stratum. These combinatorial data appear in a natural generalization of the well know relation between translation surfaces and interval exchange maps (see \cite{Masur82,Veech82}).
There is a natural construction that associates a flat surface to a generalized permutation and a continuous parameter. This is the \emph{Veech construction}.
The set of generalized permutations that can appear with this construction in a connected component of a stratum is called the \emph{extended Rauzy class}. The important fact is that it can be built canonically using the (extended) Rauzy moves. A difficulty here is that a generic surface near the boundary does not appear naturally from the Veech constrution with a ``short'' parameter (see section~\ref{ends:of:strata} for a precise statement).
We will proceed in the following way:
\begin{enumerate}
\item Given a generalized permutation $\pi$, we define a subset $Z(D_{\pi,\varepsilon})$, of flat surfaces that are near the boundary, were we take ``short'' parameters.
\item For any flat surface near the boundary, there is a path that stays near the boundary and reaches $Z(D_{\pi,\varepsilon})$, for some $\pi$ (Lemma~\ref{rectangle}).
\item We show that the subset $\cup_{\pi} Z(D_{\pi,\varepsilon})$ is connected, where the union is taken on an extended Rauzy class (Lemma~\ref{intersection:D} and Lemma~\ref{D:connected}).
\end{enumerate}

\section{Flat surfaces and moduli space}
\subsection{Generalities}
A {\it flat surface} is a real, compact, connected surface of genus $g$ equipped with a flat metric with isolated conical singularities and such that the holonomy group belongs to $\Z / 2\Z$. Here holonomy means that the parallel transport of a vector along a loop brings the vector back to itself or to its opposite. This implies that all cone angles are integer multiples of $\pi$. We also fix a choice of a parallel line field in the complement of the conical singularities. This parallel line field will be usually referred as \emph{the vertical direction}.
 Equivalently a flat surface is a triple $(S,\mathcal U,\Sigma)$ such that $S$ is a topological compact connected surface, $\Sigma$ is a finite subset of $S$ (whose elements are called {\em singularities}) and $\mathcal U = \{(U_i,z_i)\}$ is an atlas of $S \setminus \Sigma$ such that the transition maps $z_j \circ z_i^{-1} : z_i(U_i\cap U_j) \rightarrow z_j(U_i\cap U_j)$ are translations or half-turns: $z_i = \pm z_j + c$, and for each $s\in \Sigma$, there is a neighborhood of $s$ isometric to a Euclidean cone. Therefore, we get a {\it quadratic differential} defined locally in the coordinates $z_i$ by
the formula $q=d z_i^2$. This form extends to the points of $\Sigma$ to zeroes, simple poles or marked points (see~\cite{Masur:Tabachnikov,Hubbard:Masur}). 

Observe that the holonomy is trivial if and only if there exists a sub-atlas such that all transition functions are translations, or equivalently if the quadratic differential $q$ is the global square of an Abelian differential. We usually say that $S$ is a translation surface. In this case, we can choose a parallel vector field instead of a parallel line field, which is equivalent to fix a square root of $q$. In the complementary case, we sometime speak of \emph{half-translation} surfaces.

The moduli space of quadratic differentials on a Riemann surface of genus $g$ is naturally stratified by considering the quadratic differentials that have prescribed orders of zeroes (and poles). Each corresponding stratum is a complex analytic orbifold. 

A saddle connection is a geodesic segment (or geodesic loop) joining two singularities (or a singularity to itself) with no singularities in its interior. Even if $q$ is not globally a square of an Abelian differential, we can find a square root of $q$ along any saddle connection. Integrating $q$ along the saddle connection we get a complex number
(defined up to multiplication by $-1$). Considered as a planar vector, this complex number represents the affine holonomy vector along the saddle connection. In particular, its Euclidean length is the modulus of its holonomy vector. Note that a saddle connection persists under any small deformation of the surface. 
Local coordinates in a stratum of quadratic differentials are obtained by considering affine holonomy vectors of a well chosen collection of saddle connections.

A \emph{normalized stratum} corresponds to flat surfaces with area one.
There is a natural action of $\textrm{SL}_2(\mathbb{R})$ on  any normalized stratum: let $(U_i,z_i)_{i\in I}$ be an atlas of flat coordinates of $S$, with~$U_i$ open subsets of $S$ and $z_i(U_i)\subset \mathbb{R}^2$. For $A\in \textrm{SL}_2(\mathbb{R})$, an atlas of $A.S$ is given by $(U_i,A\circ \phi_i)_{i\in I}$. The action of the diagonal subgroup of $\textrm{SL}_2(\mathbb{R})$ is called the Teichmüller geodesic flow. In order to specify notations, we denote by $g_t$, $h_t$ the matrices
\begin{eqnarray*}
g_t=
\left[
\begin{array}{cc}
 e^{\frac{t}{2}} &  0    \\
 0  &    e^{-\frac{t}{2}}
\end{array} 
\right] \textrm{ and } h_t =\left[
\begin{array}{cc}
 1 &  t    \\
 0  &    1
\end{array} 
\right].
\end{eqnarray*}
It is well known (see \cite{Masur82,Veech82,Veech86}) that the Teichmüller flow and the $SL_2(\mathbb{R})$ action 
preserve a natural finite volume measure and are ergodic with respect to this measure for each connected component of each normalized stratum.

\subsection{Neighborhood of the boundary of a stratum}

Let $\CC_1$ be a connected component of a normalized stratum of the moduli space of quadratic differentials. 
Let $K_\varepsilon\subset \CC_1$ be the set corresponding to flat surfaces whose lengths of saddle connections are all bigger than or equal to $\varepsilon$. Since the set of holonomy vectors of saddle connections is discrete, we clearly have  $\cup_{\varepsilon\in \Q^+} K_\varepsilon=\CC_1$. Also, it is well known that $K_\varepsilon$ is compact.

We will call the $\varepsilon-$\emph{boundary} of $\CC_1$ the subset $\CC_{1,\varepsilon}=\CC_1 \backslash K_\varepsilon$. 

Recall that for a $\sigma$-locally compact space $W$, the number of ends is the maximal number of unbounded components of $W\backslash K$, for $K\subset W$ compact, when this number is bounded from above (see \cite{HR}). 
Hence, in order to prove the main result,  it is enough to show that the $\varepsilon$-boundary of $\CC_1$ is connected.

\section{The Rauzy--Veech induction}
In this section, we present a generalization of the Veech construction and Rauzy-Veech induction for the case of quadratic differentials. For more details, see \cite{Veech82,Marmi:Moussa:Yoccoz} for the case of Abelian differentials and \cite{B2,Boissy:Lanneau} for the case of quadratic differentials.
\subsection{Veech zippered rectangle construction}
\begin{defs}
A \emph{generalized permutation}, is a two-to-one map $\pi:\{1,\ldots,2d\}\rightarrow \{1,\ldots,d\}$, which is represented by a table of two lines of symbols, with each symbol appearing exactly two times.
$$
\pi=\left(\begin{array}{ccc}
\pi(1) & \dots& \pi(l) \\
\pi(l+1) & \dots & \pi(l+m)
\end{array}\right). 
$$
The \emph{type} of a generalized permutation is the pair $(l,m)$, where $l$ is the number of elements of the first line and $m$ is the number of element of the second line. We clearly have $l+m=2d$.

 A generalized permutation is called \emph{reduced} if for each $k$, the first occurrence in $\{1,\ldots,2d\}$ of the label $k\in\{1,\ldots,d\}$ is before the first occurrence of any label $k'>k$. 
\end{defs}

A renumbering of a generalized permutation corresponds to replacing $\pi$ by $f\circ \pi$, for $f$ a permutation of $\{1,\ldots,d\}$. In this paper, we look at  generalized permutations defined up to renumbering. It is easy to show that for each generalized permutation $\pi$, there exists a unique reduced generalized permutation $\pi_r$ which is obtained from $\pi$ by renumbering.
 
 \begin{Example}
 Consider the generalized permutation $\pi=\left(\begin{smallmatrix} 3 & 4 & 1 & 2 & 1\\ 4 & 2 & 5 & 5 & 3
\end{smallmatrix}\right)$. It is not reduced since the first occurence of the number 3 is before the first occurence of 1. In order to get a reduced generalized permutation, we clearly must replace 3 by 1, 4 by 2, etc\ldots The corresponding reduced generalized permutation therefore is $\left(\begin{smallmatrix} 1 & 2 & 3 & 4 & 3\\ 2 & 4 & 5 & 5 & 1
\end{smallmatrix}\right)$.
\end{Example}

If for all $k\leq l$ the unique $k'\neq k$ such that $\pi(k)=\pi(k')$ satisfies the condition $k'>l$, then a reduced generalized permutation satisfies $\pi(k)=k$ for all $k\leq d$, and corresponds to a permutation of $\{1,\ldots,d\}$.

\begin{defs}\label{def:gp}
Let $\pi$ be a generalized permutation of type $(l,m)$. Let $\{\zeta_{k}\}_{k \in
\{1,\ldots,d\}}$ be a collection of complex numbers such that: 
\begin{enumerate} 
\item $\forall 1\leq i \leq l-1 \quad Im(\sum_{j\leq i} \zeta_{\pi(j)})>0$ 
\item $\forall 1\leq i \leq m-1 \quad Im(\sum_{1\leq j\leq i} \zeta_{\pi(l+j)})<0$
\item $\sum_{1\leq i\leq l} \zeta_{\pi(i)} = \sum_{1\leq j\leq m}\zeta_{\pi(l+j)}$.
\end{enumerate} 
The collection $\zeta=\{\zeta_{i}\}_{i\in \{1,\ldots,d\}}$ is called a \emph{suspension datum} for $\pi$. We  will say that $\pi$ is \emph{irreducible} if $\pi$ admits suspension data.
\end{defs}

For the case when $\pi$ corresponds to a permutation, \emph{i.e.} when $l=m=d$ and $\pi(\{1,\dots ,d\})=\pi(\{d+1,\dots ,2d\})$, being irreducible means that there exists no $1\leq k< d$ such that $\pi(\{d+1,\ldots,d+k\})=\pi(\{1,\ldots k\})$. In the general case, the combinatorial criterion for being irreducible is more complicated (see \cite{Boissy:Lanneau}).

\begin{figure}[htb]
\psfrag{a}{$\scriptstyle \zeta_1$} \psfrag{b}{$\scriptstyle \zeta_2$}
\psfrag{c}{$\scriptstyle \zeta_3$} \psfrag{d}{$\scriptstyle \zeta_4$}
\psfrag{e}{$\scriptstyle \zeta_5$}

 \begin{center}
 \includegraphics[scale=0.9]{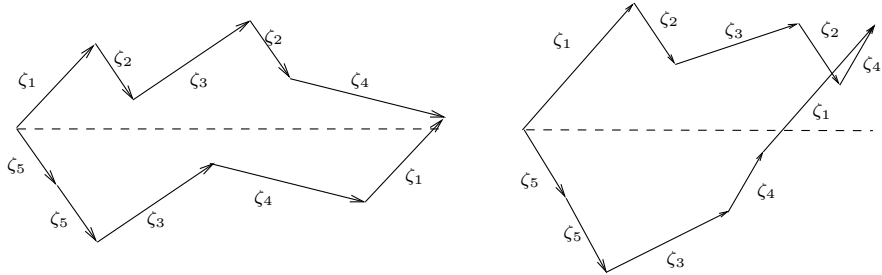}
 \caption{Examples of suspension data associated to a generalized permutation}
 \label{polygone}
 \end{center}
\end{figure}

Now we will associate to any pair $(\pi,\zeta)$ a flat surface. 
This construction will be refered as the \emph{Veech construction}\footnote{Note that the original paper of Veech \cite{Veech82} does not explicitely use suspension datum as parameter. This point of view is due to Marmi, Moussa and Yoccoz  \cite{Marmi:Moussa:Yoccoz} and is equivalent to Veech's original construction.}.
 We first describe a simple case.
We consider a broken line $L_1$ whose edge number $i$ ($1\leq i \leq l$) is represented  by the complex number $\zeta_{\pi(i)}$. Then we consider a second broken line $L_2$ which starts from the same point, and whose edge number $j$ ($1 \leq j\leq m$) is represented by $\zeta_{\pi(l+j)}$. The last condition of Definition~\ref{def:gp} implies that these two lines also end at the same point. 
If they have no other intersection points, then they form a polygon (see Figure~\ref{polygone}). The sides of the polygon, enumerated by indices of the corresponding complex numbers, come naturally as pairs of segments that are  parallel and of same length. Gluing these pairs of sides by isometries respecting the natural orientation of the polygon, this construction defines a flat surface. The holonomy of this surface is either trivial or $\mathbb{Z}/2\mathbb{Z}$, depending on the generalized permutation $\pi$.

 \begin{figure}[htb]
 \begin{center}
\psfrag{r1}{$R_1$}  \psfrag{r2}{$R_2$}  \psfrag{r3}{$R_3$}   \psfrag{r4}{$R_4$}
 \includegraphics[scale=0.6]{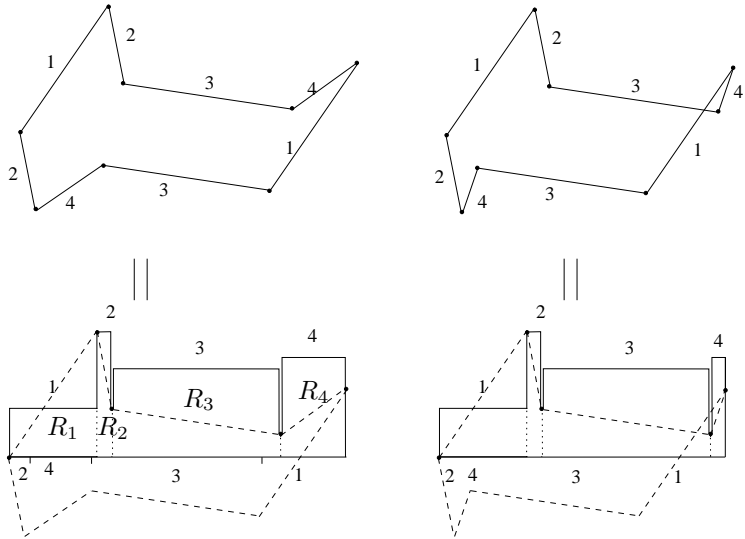}
 \caption{Veech zippered rectangle construction, for  two examples of  suspension data, in the case of Abelian differentials}
 \label{polygotozipp}
 \end{center}
\end{figure}

The two lines $L_1$ and $L_2$ might intersect in some other point, even if we start from a suspension data (see Figure~\ref{polygone}). However, we can define an alternative construction that associates, to any suspension data, a flat surface. This construction coincides with the first one when the lines $L_i$ have no other intersecting points. 
 The general idea is to consider a collection of (oriented) rectangles $(R_\alpha)_{\alpha\in \{1,\dots ,d\}}$, of width $Re(\zeta_\alpha)$ and height $h_\alpha$, where $h_\alpha>0$ depends linearly on $(Im(\zeta_\beta))_{\beta\in \{1,\dots ,d\}}$. The surface is obtained from $\sqcup_\alpha R_\alpha$ by identifications of the boundaries of the $R_\alpha$. The identifications of the horizontal boundary are done in the following way:
\begin{itemize}
\item We consider the oriented interval $I=[0,\sum_{i=1}^l Re(\zeta_{\pi(i)})]$.
\item For each $R_\alpha$, we choose a labeling of  each horizontal sides by the preimages of $\alpha$ by $\pi$.  
\item For each $i\in \{1,\dots ,l\}$, we identify with an isometry the horizontal side of $R_{\pi(i)}$, labeled by $i$,  with the segment $I_i=[\sum_{j<i} Re(\zeta_{\pi(j)}),\sum_{j\leq i} Re(\zeta_{\pi(j)})]\subset I$ as in Figure~\ref{polygotozipp}. This gluing is done so that the orientation of $I$ and the one of $\partial R_{\pi(i)}$ coincide.    
\item For each $i\in \{l+1,\dots ,l+m\}$, we identify as before with an isometry an horizontal side of $R_{\pi(i)}$ with $I_i=[\sum_{l+1\leq j<i} Re(\zeta_{\pi(j)}),\sum_{l+1\leq j\leq i} Re(\zeta_{\pi(j)})]\subset I$ as in Figure~\ref{polygotozipp}. This gluing is done so that the orientation of $I$ and the one of $\partial R_{\pi(i)}$ are opposite.
\end{itemize}
The identifications of the vertical boundary are done  in the following way:
\begin{itemize}
\item For each $i\in \{1,\dots ,l-1\}$ we ``zip'' the adjacent rectangles $R_{\pi(i)}$ and $R_{\pi(i+1)}$ starting from their common point on $I$, on a distance of $\sum_{j\leq i} Im(\zeta_{\pi(j)})$. 
\item For each $i\in \{l+1,\dots ,l+m-1\}$ we ``zip'' the adjacent rectangles $R_{\pi(i)}$ and $R_{\pi(i+1)}$ starting from $I$ on a distance of $-\sum_{l+1\leq j\leq i} Im(\zeta_{\pi(j)})$. 
\end{itemize}
 The $h_k$ are taken so that we obtain a closed surface. One can show that suitable $h_k$ exist if and only if and only if there exists a suspension data.

 On Figure~\ref{polygotozipp}, we can see the details of identifications for two examples of suspension data.
This more general construction is called the \emph{Veech zippered rectangle construction}. We refer to \cite{Veech82,Marmi:Moussa:Yoccoz} for the case of Abelian differentials and to \cite{B2,Boissy:Lanneau} otherwise for more details.

In all cases, we will denote by $Z(\pi,\zeta)$ the flat surface obtained from this construction. We will need the following lemma:

\begin{lem}\label{scmin}
Let $\pi$ be a irreducible generalized permutation and let $\zeta$ be a suspension datum for $\pi$. The flat surface $Z(\pi,\zeta)$ obtained by the Veech construction has a saddle connection of length less than or equal to $\min\{|\zeta_\alpha|, \alpha\in \mathcal{A}\}$.
\end{lem}

\begin{proof}
In the case when the two lines $L_1$ and $L_2$ intersect only in their endpoints, it is clear that the segments of the lines $L_1,L_2$ will corresponds to saddle connections. Hence, each $|\zeta_\alpha|$ is the length of a saddle connection for each $\alpha$ and the lemma is proven.

Now, let us assume that the two lines intersect elsewhere. Then we must have $\sum_{k=1}^l Im(\zeta_{\pi(k)})\neq 0$ otherwise there would be no other intersection points. Without loss of generality, one can assume that this sum is positive. The other intersection points correspond to the intersection of $L_1$ with  the rightmost segment of $L_2$, since it is the only one which is not a subset of $\{z\in \mathbb{C}, Im(z)\leq 0\}$. In particular, $\zeta_{\pi(l+m)}$ does not correspond to a saddle connection. However, all the other parameters $\zeta_\alpha$, for $\alpha\neq \pi(l+m)$ correspond to a saddle connection which is inside the corresponding rectangle in the zippered rectangle construction. 

\begin{figure}[htb]
\begin{center}
\psfrag{g}{$\gamma$} 
\psfrag{z}{$\zeta_{\pi(l+m)}$}

 \includegraphics[scale=0.6]{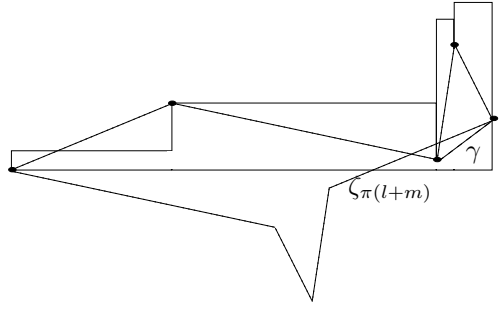}
 \caption{The dashed segment $\gamma$ corresponds to a saddle connection of lenght less than $|\zeta_{\pi(l+m)}|$}
 \label{badpolysc}
 \end{center}
\end{figure}
  
\end{proof}

 Let $k_0$ be the maximal index smaller than $l$ such that the vertex of $L_1$ of coordinates $\sum_{i\leq k_0} \zeta_{\pi(i)}$ is under the line $L_2$. Let $\gamma$ be the segment in $\mathbb{C}$ that joins the two points of coordinates $\sum_{i\leq k_0} \zeta_{\pi(i)}$ to $\sum_{i\leq l} \zeta_{\pi(i)}$. By hypothesis on $k_0$, for each $k_0<k<l$, the point in the plane of coordinates $\sum_{i\leq k} \zeta_{\pi(i)}$ is above the point of $\gamma$ with the same real part. Hence, it is easy to see that $\gamma$ corresponds to a saddle connection on the surface $Z(\pi,\zeta)$ (see Figure~\ref{badpolysc}). Also, $\gamma$ is of length smaller than $\zeta_{\pi(l+m)}$ since its real part and imaginary part are necessarily smaller than the ones of $\zeta_{\pi(l+m)}$.

\subsection{Rauzy--Veech induction and extended Rauzy classes}
Now we will define the Rauzy--Veech induction on the space of suspension data. 
Let us fix some terminology: let $k\in \{1,\ldots,l+m\}$. The \emph{other occurrence} of the symbol $\pi(k)$ is the unique integer $k'\in \{1,\ldots,l+m\}$, distinct from $k$, such that $\pi(k')=\pi(k)$. 

Now we define the combinatorial Rauzy moves for the generalized permutations. We first define the unreduced maps $\mathcal{R}_0'$, $\mathcal{R}_1'$ and $s'$.
\begin{enumerate}
\item $\mathcal R_0'(\pi)$: \\
$\bullet$ If the other occurrence $k$ of the symbol $\pi(l)$ is in $\{l+1,\ldots,l+m-1\}$, then we define $\mathcal R'_0(\pi)$ to be of type $(l,m)$ obtained by removing the symbol $\pi(l+m)$ from the occurrence $l+m$ and putting it at the occurrence $k+1$, between the symbols $\pi(k)$ and $\pi(k+1)$. \\
$\bullet$ If the other occurrence $k$ of the symbol $\pi(l)$ is in $\{1,\ldots,l-1\}$, and if there exists another symbol $\alpha$, whose both occurrences are in $\{l+1,\ldots,l+m\}$, then we we define $\mathcal R'_0(\pi)$ to be of type $(l+1,m-1)$ obtained by removing the symbol $\pi(l+m)$ from the occurrence $l+m$ and putting it at the occurrence $k$, between the symbols $\pi(k-1)$ and $\pi(k)$ (if $k=1$,  by convention the symbol $\pi(l+m)$ is put on the left of the first symbol $\pi(1)$).\\
$\bullet$ Otherwise $\mathcal R'_0 \pi$ is not defined. 

\item The map $\mathcal{R}_1'$ is obtained by conjugating $\mathcal{R}_0'$ with the transformation that interchanges the two lines in the table representation.

\item The map $s'(\pi)$ is the generalized permutation of type $(m,l)$ such that  $s'(\pi)(k)=\pi(2d-k+1)$.
\end{enumerate}
Then we obtain $\mathcal{R}_0(\pi)$, $\mathcal{R}_1(\pi)$ and $s(\pi)$ by renumbering $\mathcal{R}_0'(\pi)$, $\mathcal{R}_1'(\pi)$ and $s'(\pi)$ in order to get reduced generalized permutations.
For a more explicit definition of $\mathcal{R}'_0$ and $\mathcal{R}'_1$  in terms of the map $\pi$, we refer to \cite{Boissy:Lanneau}.

\begin{Example}
Let $\pi$ be the generalized permutation 
$\pi=\left(\begin{smallmatrix} 1 & 2 & 3 & 4 & 3\\ 2 & 4 & 5 & 5 & 1
\end{smallmatrix}\right)$. We have
$$
\mathcal \mathcal{R}'_0(\pi)=\left(\begin{array}{cccccc} 1 & 2 & 1 & 3 & 4 & 3\\ 2 & 4 & 5 & 5 &&
\end{array}\right) =
\mathcal \mathcal{R}_0(\pi),
$$
 and
 $$
 \mathcal \mathcal{R}'_1(\pi)=\left(\begin{array}{cccccc} 1 & 3&2&3&4\\ 2&4&5&5&1
\end{array}\right) \textrm{ so }
\mathcal \mathcal{R}_1(\pi)=\left(\begin{array}{cccccc} 1 &2 & 3 & 2 & 4 \\ 3 & 4 & 5 & 5 & 1
 \end{array}\right).
$$

Also,  
$$s'(\pi)=\left(\begin{array}{ccccc}1&5&5&4&2 \\ 3&4&3&2&1
\end{array}\right) \textrm{ so }
s(\pi)=\left(\begin{array}{ccccc}
1&2&2&3&4\\ 5&3&5&4&1
\end{array}\right). 
$$

\end{Example}

We now define the \emph{Rauzy--Veech induction} on the space of suspension data.
\begin{itemize}
\item 
Let $\zeta$ be a suspension datum for $\pi$ such that we have $Re(\zeta_{\pi(l)})>Re(\zeta_{\pi(l+m)})$. We can define a suspension datum $\zeta'$ for $\pi'=\mathcal{R}'_0(\pi)$ in the following way:
$\zeta^{\prime}_k=\zeta_k$ if $k\neq \pi(l)$ and
$\zeta_{\pi(l)}^{\prime}=\zeta_{\pi(l)}-\zeta_{\pi(l+m)}$. 
\item
Let $\zeta$ be a suspension datum for $\pi$ such that $Re(\zeta_\pi(l))<Re(\zeta_\pi(l+m))$. We can define a suspension datum $\zeta'$ for $\pi'=\mathcal{R}'_1(\pi)$ in the following way:
$\zeta^{\prime}_k=\zeta_k$ if $k\neq \pi(l+m)$ and
$\zeta_{\pi(l+m)}^{\prime}=\zeta_{\pi(l+m)}-\zeta_{\pi(l)}$. 
\end{itemize}

\begin{figure}[htbp]
 \begin{center}
 \psfrag{1}[][]{\small $\zeta_1$}
 \psfrag{2}[][]{\small $\zeta_2$}
 \psfrag{3}[][]{\small $\zeta_3$}
 \psfrag{4}[][]{\small $\zeta_4$}

 \psfrag{1p}[][]{\small $\zeta'_1$}
 \psfrag{2p}[][]{\small $\zeta'_2$}
 \psfrag{3p}[][]{\small $\zeta'_3$}
 \psfrag{4p}[][]{\small $\zeta'_4$}

 \includegraphics[width=360pt]{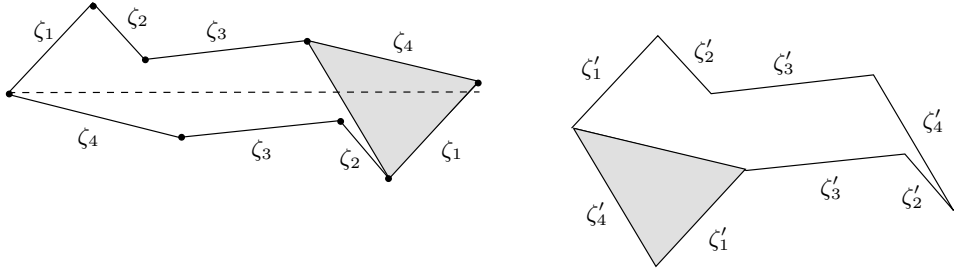}
 \caption{Rauzy-Veech induction for a suspension datum.}
 \label{fig:rauzy:suspension}
 \end{center}
\end{figure}

\begin{remark}\label{rauzy:susp}
The parameter  $\zeta^{\prime}$ defines a suspension datum for $\pi'$. The flat surfaces $Z(\pi,\zeta)$ and $Z(\pi',\zeta')$ are naturally identified in the moduli space (see Figure~\ref {fig:rauzy:suspension}).
\end{remark}

A consequence of this remark is that if $\pi$ is irreducible, then so are $\mathcal{R}_0(\pi)$ and $\mathcal{R}_1(\pi)$, and they correspond to the same connected component of the moduli space of quadratic differentials.

Note that when $\pi$ does not correspond to a permutation, $s(\pi)$ is not necessary irreducible as shown in the next example. When $s(\pi)$ is irreducible, then it corresponds to the same connected component of the moduli space of quadratic differentials.  


 \begin{Example}
 Let $\pi=\left(\begin{smallmatrix} 1& 2& 1\\2&3&3&4&4 \end{smallmatrix}\right)$. This generalized permutation is irreducible: for instance $\zeta=(2+2i,1-i,1-2i,1+4i)$ is a suspension data for $\pi$. However $s'(\pi)=\left(\begin{smallmatrix}4&4&3&3&2\\ 1&2&1\end{smallmatrix}\right)$ does not admit a suspension data. Indeed, if $\zeta'$ were a suspension datum for $s'(\pi)$, then we would get:
 $$
 Im(\zeta'_1)<0 \quad \textrm{and} \quad 2Im(\zeta'_4)+2Im(\zeta'_3)>0$$
 However, we must have $\zeta'_1=\zeta'_3+\zeta'_4$ which contradicts the previous inequalities.
 \end{Example}

\begin{defs}
\begin{enumerate}
\item A \emph{Rauzy class} is a minimal subset of irreducible generalized permutations (or permutations) which is invariant by the combinatorial Rauzy maps $\mathcal{R}_0,\mathcal{R}_1$. 
\item An \emph{extended Rauzy class} is a minimal subset of irreducible generalized permutations (or permutations) which is invariant by the combinatorial Rauzy maps $\mathcal{R}_0,\mathcal{R}_1$ and $s$.
\end{enumerate}
\end{defs}

Extended Rauzy classes can be build by the following way: we start from a generalized permutation and consider all its descendant by the Rauzy moves, were we forbid the operation $s$ on $\pi$ when $s(\pi)$ is reducible.

One have the following result, due to Veech (see \cite{Veech90}) for the case of permutations and to the author and E. Lanneau for the case of generalized permutation (see Appendix~B in \cite{Boissy:Lanneau}).

\begin{thm}\label{extrc:cc}
Extended Rauzy classes are in one to one correspondence with the connected components of the moduli space of quadratic differentials.
\end{thm}

\begin{remark}\label{rem:erc}
A sufficient condition for $s(\pi)$ to be irreducible is that $\pi$ admits a suspension data such that $$\sum_{k=1}^{l} Im(\zeta_{\pi(k)})=\sum_{k=l+1}^{l+m}Im(\zeta_{\pi(k)})=0.$$
Indeed, in this case, $\zeta$ is a suspension data for $s(\pi)$ after suitable renumbering. In the proof of Theorem~\ref{extrc:cc} in \cite{Boissy:Lanneau}, the operation $s$ is only used on such generalized permutations. Hence we can weaken the definition of extended Rauzy class by authorizing the map $s$ only for such permutations. 
\end{remark}

\section{Ends of strata} \label{ends:of:strata}
In order to prove Theorem~\ref{end}, we define a subset of $\CC_{1,\varepsilon}$ that corresponds to flat surfaces obtained from the zippered rectangle construction, and with a ``small'' parameter. We show that any connected component of $\CC_{1,\varepsilon}$ must intersect this set, and that it is connected.

The following definition is needed for technical reasons when the corresponding stratum corresponds to quadratic differentials.
\begin{defs}\label{regular}
Let $\pi$ be a generalized permutation. A symbol $\alpha$ is said to be \emph{regular} if one of the following property holds:
\begin{itemize}
\item 
 $\alpha$ appears in both lines,
 \item
  there exists $\beta\neq \alpha$ that appears only in the same line as $\alpha$.
  \end{itemize}
\end{defs}
 
The geometric interpretation of this notion is the following: for a regular symbol $\alpha$, there exists a suspension datum $\zeta$ such that $Re(\zeta_\alpha)$ is smaller than $Re(\zeta_\beta)$ for all $\beta\neq \alpha$. For a nonregular symbol, such suspension does not exists since the suspension data condition implies that $\zeta_\alpha=\sum_{\beta} \zeta_\beta$, were the sum is taken on all $\beta$ that appear only in the line different from $\alpha$.

\begin{Example}
Let $\pi$ be the generalized permutation $\pi=\left(\begin{smallmatrix} 1 & 2 & 3 & 4 & 3\\ 2 & 4 & 5 & 5 & 1 &6 &6
\end{smallmatrix}\right)$.
 The symbols ''1'', ''2'' and ''4'' satisfy the first condition, and the symbols ``5'', ``6'' satisfy the second condition, hence they are regular. The symbol ''3'' satisfies none of the prescibed condition, hence it is not regular.
\end{Example}

\begin{defs}\label{def:D_pi}
Let $\pi$ be  a generalized permutation.  We denote by  $D_{\pi,\varepsilon}$ the set of of pairs $(\pi,\zeta)$ , with $\zeta$ a suspension data for $\pi$, that define an area one surface and such that there exists a regular symbol $\alpha\in \{1,\ldots,d\}$ with $|\zeta_\alpha|< \varepsilon$.
\end{defs}

Recall that we denote by $Z(\pi,\zeta)$ the flat surface obtained from $(\pi,\zeta)$ by the Veech construction. By Lemma~\ref{scmin}, we have $Z(D_{\pi,\varepsilon})\subset \CC_{1,\varepsilon}$. Conversely, given a surface in $\CC_{1,\varepsilon}$, one would like to present it as in some $Z(D_{\pi,\varepsilon})$. However, it is not necessarily the case: even when the surface is obtained by the Veech construction, there is not necessarily a small parameter in the corresponding suspension data. So we will use the following lemma:

\begin{lem}\label{rectangle}
Let $S\in \CC_{1,\varepsilon}$. We assume that the underlying stratum is neither  $\mathcal{H}(\emptyset)$, nor $\mathcal{Q}(-1,-1,-1,-1)$.  There exists a generalized permutation $\pi$ and a suspension datum $\zeta$ such that $(\pi,\zeta)\in D_{\pi,\varepsilon}$ and 
 $Z(\pi,\zeta)$ is in the same connected component of $\CC_{1,\varepsilon}$ as $S$.
\end{lem}

\begin{proof}
We claim that we can find a surface $S'$ in the same connected component of $\mathcal{C}_{1,\varepsilon}$ as $S$,  whose horizontal foliation consists of one cylinder, and with a horizontal saddle connection of length smaller than~$\varepsilon$.

The set of flat surfaces with a horizontal one cylinder decomposition is dense in each connected component of a stratum (see \cite{Kontsevich:Zorich}, Remark~7 and \cite{Lanneau:cc}, Theorem~3.6). Hence, up to a small perturbation of $S$ we can assume that it is the case for  $S$. It admits a smallest saddle connection $\gamma$. If $\gamma$ is horizontal, then the claim is true.  We  assume now that $\gamma$ is not horizontal. Then, we apply on the surface the matrix $h_t=\left( \begin{smallmatrix} 1 &t\\ 0 &1 \end{smallmatrix} \right)$ for a suitable continuous family of parameter $t$, so that $\gamma$ becomes vertical. Note that in the last operation, the area of the flat surface does not change and the size of $\gamma$ decrease, and therefore the flat surface stays in $\CC_{1,\varepsilon}$ during that process.

The resulting surface, that we still denote by $S$ can be represented by a rectangle whose vertical sides are identified, and  whose horizontal sides admit a partition in pairs of segments, and for each pair, the segments are of the same length and are identified either by translation or by a half-turn (see Figure~\ref{longcyl}). The endpoints of each segments correspond to singularities of $S$. The bottom left corner can be assumed to be a singularity which is an end point of $\gamma$, but the top left point does not correspond in general to a singularity. Indeed, for a generic initial surface, the corresponding rectangle is very long and very thin, and the saddle connection $\gamma$ is much bigger than the heigth of that rectangle. In particular, in the notation of Figure~\ref{longcyl},  the segments labeled ``$2a$'', and ``$2b$'' have a nonsingular endpoint, should be seen as the two parts of a single a horizontal segment. Let $N$ be ratio of the length of $\gamma$ by the heigth $\eta$ of the rectangle. This ratio $N$ is necessarily an integer.

\begin{figure}[htb]
\begin{center}
\psfrag{0}{\tiny 0}
\psfrag{1}{\tiny 1}
\psfrag{2a}{\tiny 2a}
\psfrag{2b}{\tiny 2b}
\psfrag{3}{\tiny 3}
\psfrag{4}{\tiny 4}
\psfrag{5}{\tiny 5}
\includegraphics[scale=0.5]{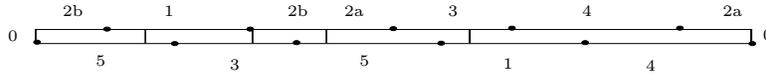}
\caption{A small vertical saddle connection on a long thin cylinder}
\label{longcyl}
\end{center}
\end{figure}

Now we assume that there is a pair $(s_1,s_2)$ of horizontal segments on the rectangle  that are identified by a translation. Note that $s_1$ or $s_2$ might be the concatenation of the top-left and top-right horizontal segments 
We want to shrink $s_1$ and $s_2$ until they are very small. During the process, we must modify the surface so that it stays in $\CC_{1,\varepsilon}$.

Let $A<1$ be the area of the limit rectangle obtained after collapsing the segments $s_1,s_2$ to a point. Note that since $S\notin \mathcal{H}(\emptyset)$ by asumption, there exists other horizontal segments than $s_1,s_2$ on the initial rectangle, and therefore $A\neq 0$. Up to applying the matrix $g_{t_0}$ for a suitable $t_0$ to the surface $S$, we can assume that  $N.\eta<\sqrt{A}\varepsilon$. Now we continuously shrink the segments $s_1$ and $s_2$ until they are small enough.  We get a  continuous family of flat surfaces $(\tilde{S}_t)$. There exists $t_1\in \mathbb{R}$ such that the saddle connection $\gamma$ persists for $t\in [0,t_1[$, but does not necessarily stay vertical.  We can find a continuous function $\varphi$ such that the saddle connection corresponding to $\gamma$ in $h_{\varphi(t)}.\tilde{S}_t$ stays vertical. Rescaling $h_{\varphi(t)}.\tilde{S}_t$ so that it has area one, we get a continuous family of surfaces $(S_t)_t$ in $\CC_{1,\varepsilon}$  with a one cylinder decomposition, and with a vertical saddle connection of length $N.\eta_t$, where $\eta_t$ is the heigth of the rectangle. Since the area $A_t$ of $h_{\varphi(t)}.\tilde{S}_t$ is bigger than $A$, and the rescaling changes the lengths of all saddle connections by a factor $1/\sqrt{A_t}$, then we have $N.\eta_t<\varepsilon$. Hence $S_t$ is in $\CC_{1,\varepsilon}$ for all $t< t_1$.

Taking $t_1$ maximal, either we find $t<t_1$ such that the lenght of $s_i$ is smaller than $\varepsilon$, and therefore the corresponding surface $S_t$ has a horizontal saddle connection of length smaller than $\varepsilon$, or $S_t$ converges to a surface $S_{t_1}$ when  $t$ tends to $t_1$. The path corresponding to $\gamma$ in $S_{t_1}$ is not a saddle connection but a union of saddle connections. This means that we have a saddle connection $\gamma'_{t_1}$, of length $N'\eta_{t_1}<\varepsilon$ with $N'$ an integer smaller than $N$. And we can continue the shrinking process. There can be only a finite number of such steps, hence after a while, the segments $s_1,s_2$ are of length smaller than $\varepsilon$, thus we have found the surface $S'$ of the claim.

If the pair $(s_1,s_2)$ does not exist then the surface is a half-translation surface and there exists two pairs of segments $(s_1,s_2)$ and $(s_1',s_2')$, one on each horizontal side of the rectangle, whose corresponding identifications are half-turns. The previous proof works if we shrink both of them at the same speed until the smallest one, say $(s_1,s_2)$, becomes small enough. Here the hypothesis that $S\notin \mathcal{Q}(-1,-1,-1,-1)$ implies that there is at least another pair of horizontal segment in the rectangle, on the same side of the rectangle as the pair $(s_1,s_2)$. In particular, the limit area  $A$ is nonzero, since the limit rectangle is nontrivial. 

Now, starting from $S'$, we can continuously change the vertical direction and we get a area one surface obtained by the zippered rectangle construction with data $(\pi,\zeta)$  such that $|\zeta_\alpha|<\varepsilon$ for some $\alpha$.
In the first case of the proof $\alpha$ is clearly regular. In the second case, the small $\zeta_\alpha$ corresponds to the smallest pair of horizontal segments $(s_1,s_2)$. By the above remark, there is another pair of horizontal segments on the same side as $(s_1,s_2)$ that are identified together by an half-turn. This precisely means that $\alpha$ is regular.
\end{proof}

Now we want to prove that $\cup_{\pi} Z(D_{\pi,\varepsilon})$ is contained in the same connected connected component of $\mathcal{C}_{1,\varepsilon}$, where the union is taken over all the elements in the extended Rauzy class. This is a consequence of the next two lemmas.

\begin{lem}\label{intersection:D}
 Let $\pi$ be an irreducible permutation or generalized permutation. 
 \begin{enumerate}
\item Let $\pi'$ be $\mathcal{R}_0(\pi)$ or $\mathcal{R}_1(\pi)$.  Then $Z(D_{\pi,\varepsilon})\cap Z(D_{\pi',\varepsilon})\neq \emptyset$.
\item If $\pi$ admits a suspension datum such that $\sum_{k\leq l} Im(\zeta_{\pi(k)})=0$ then $Z(D_{\pi,\varepsilon})$ and $Z(D_{s(\pi),\varepsilon})$ have elements in the same connected component of $\mathcal{C}_{1,\varepsilon}$.
\end{enumerate}
\end{lem}

\begin{proof}
Let $(l,m)$ be the type of $\pi$. 

$1)$ We assume that the step of the Rauzy induction we consider is $\mathcal{R}_1$. This means in particular that $\pi(l)$ is regular.
We first build $(\pi,\zeta)\in D_{\pi,\varepsilon}$ such that $|\zeta_{\pi(l)}|<\varepsilon$ and $Re(\zeta_{\pi(l)})<Re(\zeta_{\pi(l+m)})$. 

Since $\pi$ is irreducible and $\pi(l)$ is regular, there exists a suspension data $\zeta$ such that $Re(\zeta_{\pi(l)})$ is smaller than $Re(\zeta_{\pi(l+m)})$. 
Now we want to shrink $Re(\zeta_{\pi(l)})$ so that, after area renormalization, we can assume that the product  $Re(\zeta_{\pi(l)}).Im(\zeta_{\pi(l)})$ is smaller than $\varepsilon^2 /4$. If the occurrences of $\pi(l)$ appear on both sides, this is obviously possible. If the occurrences of $\pi(l)$ appear only on the top side, we must fulfill the condition $\sum_{i\leq l}Re(\zeta_{\pi(i)})=\sum_{i>l} Re(\zeta_{\pi(i)})$. This is done by increasing  
 $Re(\zeta_\beta)$ for some other (possibly all) symbols  that appear only on the top side. Such symbols exists because $\alpha$ is regular.
 
 Once the condition $Re(\zeta_{\pi(l)}).Im(\zeta_{\pi(l)})<\varepsilon^2 /4$ is obtained,  we can apply the Teichmüller geodesic flow $g_t$ on the corresponding flat surface until we have $Re(\zeta_{\pi(l)})<\varepsilon/2$ and $Im(\zeta_{\pi(l)})<\varepsilon/2$. Then $|\zeta_{\pi(l)}|<\varepsilon$ and $Re(\zeta_{\pi(l+m)})>Re(\zeta_{\pi(l)})$. Since $\pi(l)$ is regular, $(\pi,\zeta)\in D_{\pi,\varepsilon}$.

Applying the Rauzy--Veech induction on $(\pi,\zeta)$, we obtain $(\pi',\zeta')$ that define the same surface $S$ and there is a symbol $\alpha$ such that $\zeta'_\alpha=\zeta_{\pi(l)} $. It is easy to check that this symbol is regular for $\pi'$, and therefore $(\pi',\zeta')\in D_{\pi',\varepsilon}$ and $Z(\pi,\zeta)=Z(\pi',\zeta')\in Z(D_{\pi,\varepsilon})\cap Z(D_{\pi',\varepsilon})$. This proves the first point of the lemma for the map $\mathcal{R}_1$. The other case is symmetric.
\medskip

$2)$
Now we start from a suspension datum with $\sum_{k\leq l} Im(\zeta_{\pi(k)})=0$.  First, we remark that the existence of a suspension data implies that either the symbol $\pi(l)$ or the symbol $\pi(l+m)$ is regular. Indeed, otherwise the suspension data condition $\sum_{k=1}^l \zeta_{\pi(k)}=\sum_{k=l+1}^{l+m}\zeta_{\pi(k)}$ implies $2\zeta_{\pi(l)}=2\zeta_{\pi(l+m)}$, and therefore $\sum_{k=1}^{l-1} Im(\zeta_{\pi(k)})=\sum_{k=l+1}^{l+m-1}Im(\zeta_{\pi(k)})$ is both positive and negative, which is impossible. Therefore, at least one of the Rauzy moves $\mathcal{R}_0$, $\mathcal{R}_1$ is always possible. Then the same construction as before gives $(\pi,\zeta')\in D_{\pi,\varepsilon}$ such that $\sum_{k\leq l} Im(\zeta_{\pi(k)}')=0$. 
We claim that the parameter $\zeta'$ is also suspension datum for the unreduced generalized permutation $s'(\pi)$. Indeed, by definition, $s'(\pi)(k)=l+m-k+1$ and $s'(\pi)$ is of type $(m,l)$. Then, for $k< m$
\begin{eqnarray*}
\sum_{i=1}^k Im(\zeta_{s'(\pi)(i)})&=&\sum_{i=l+m-k+1}^{l+m} Im(\zeta_{\pi(i)})\\
&=&\sum_{i=l+m-k+1}^{l+m} Im(\zeta_{\pi(i)})-\sum_{i=l+1}^{l+m} Im(\zeta_{\pi(i)})\\
&=& -\sum_{i=l+1}^{l+m-k} Im(\zeta_{\pi(i)})>0
\end{eqnarray*}
Similarly, we have $\sum_{i=m+1}^{m+k}Im(\zeta_{s'(\pi)(i)})<0$ for all $1\leq k<l$.

If the underlying stratum is a not a stratum of Abelian differentials, then we have $Z(\pi,\zeta')=Z(s'(\pi),\zeta')$, and after renumbering, we have $Z(\pi,\zeta)\in Z(D_{s(\pi),\varepsilon})$. If the underlying stratum is a stratum of Abelian differentials, then we consider the path $\left(r_{\theta}.Z(\pi,\zeta)\right)_{\theta\in [0,\pi]}$ for $r_\theta=\left(\begin{smallmatrix} \cos(\theta)& -\sin(\theta)\\\sin(\theta)& \cos(\theta)
\end{smallmatrix}\right)$. The latter path is clearly in $\mathcal{C}_{1,\varepsilon}$ and joins $Z(\pi,\zeta')$ and $Z(s'(\pi),\zeta')$. This proves the lemma.

\end{proof}

\begin{lem}\label{D:connected}
Let $\pi$ be a generalized permutation whose corresponding stratum is neither  $\mathcal{H}(\emptyset)$, nor $\mathcal{Q}(-1,-1,-1,-1)$. The set $D_{\pi,\varepsilon}$ is connected.
\end{lem}

\begin{proof}
We consider $\zeta,\zeta'$ that define two elements in $D_{\pi,\varepsilon}$. There exists $\alpha,\alpha'$ in $\{1,\ldots, d\}$ regular such that $|\zeta_\alpha|<\varepsilon$ and  $|\zeta'_{\alpha'}|<\varepsilon$ (we may have $\alpha=\alpha'$).

 \emph{Case 1}: We assume that there is a symbol $\beta$, with $\beta\neq \alpha$ and $\beta \neq \alpha'$ that appears in both lines of $\pi$.

Let $\zeta_t$ be the suspension data obtained after rescaling $(1-t) \zeta+t \zeta'$ so that the corresponding flat surface has area one. We denote by $h_t$ height of the rectangle  associated to the  label $\beta$. Since $t\mapsto \zeta_t$ is continuous on $[0,1]$, there exists $h>0$ such that $h_t>h$ for all $t$ and there exists $c>0$ such that  $|\zeta_{t,\alpha}|<c$ and $|\zeta_{t,\alpha'}|<c$ for all $t$.

Let $N>0$. We consider $\zeta^{(N)}_{t}$ obtained from $\zeta_t$ after adding $N$ to the parameter $\zeta_{t,\beta}$. This means that for each $t$ we increase the length of the intervals $I_k,I_{k'}$ of the corresponding linear involution by $N$. By construction, the area of the flat surface defined by $\zeta^{(N)}_{t}$ is $A_t=1+Nh_t$. 

Then,  $\tilde{\zeta}_t^{(N)}=\frac{1}{\sqrt{A_t}}.\zeta^{(N)}_{t}$ is  a suspension data and we have, for all $t$:
$$
|\tilde{\zeta}_{t,\alpha}^{(N)}|< \frac{c}{\sqrt{1+Nh_t}}< \frac{c}{\sqrt{1+Nh}}
$$
since $h_t>h$, and similarly
$$|\tilde{\zeta}_{t,\alpha'}^{(N)}|< \frac{c}{\sqrt{1+Nh}}$$
For $N$ large enough, $\frac{c}{\sqrt{1+Nh}}$ is smaller than $\varepsilon$ and hence $(\pi,\tilde{\zeta}_t)$ is in $D_{\pi,\varepsilon}$ for all $t\in [0,1]$.

Furthermore, rescaling suitably the path $\left(\zeta_{0}^{(sN)}\right)_{s\in [0,1]}$ (respectively, the path $\left(\zeta_1^{(sN)}\right)_{s\in [0,1]})$, we see that  $(\pi,\zeta)$ and $(\pi,\tilde{\zeta})$ (resp. $(\pi,\zeta')$ and $(\pi,\tilde{\zeta'})$) are in the same connected component of $D_{\pi,\varepsilon}$. This proves the lemma for Case~1.
\medskip

\emph{Case 2:}  We assume that there are two symbols $\beta,\beta'$ different from $\alpha,\alpha'$, and such that $\beta$ appears only in the top line of $\pi$, and $\beta'$ appears only in the bottom line. The same proof as before works if we add $N$ to $\zeta_{t,\beta}$ and $\zeta_{t,\beta'}$.

\medskip
Now we come back  to the general case.  For Abelian differentials, the Case~1 is always true except if $d=2$, but this corresponds to the stratum $\HH(\emptyset)$. For quadratic differentials, recall that $\alpha$ and $\alpha'$ are regular. There are several possibilities:
\begin{itemize}
\item Each of the symbols $\alpha$ and $\alpha'$ appears in both lines of $\pi$, then Case~2 is satisfied.
\item The symbol $\alpha$ appears only in one line (for instance the top one) and the symbol $\alpha'$ appears in both.  Then, there is $\beta\neq \alpha$ that appears only in the top line because $\alpha$ is regular. Since there exists necessarily $\beta'$ that appears only in the second line, Case~2 applies.
\item  The symbol $\alpha$ appears only in the top line and the symbol $\alpha'$ appears only in the bottom line. Then, Case~2 applies since $\alpha$ and $\alpha'$ are regular.
\item Both symbols $\alpha,\alpha'$ appear only in one line (for instance the top one). We can assume that none of the previous assumption is valid. In particular there are no other symbols in the first line. If there is only one symbol in the second line, then the generalized permutation is defined with only three symbols and therefore one can check that the stratum is $\QQ(-1,-1,-1,-1)$, contradicting the hypothesis.
 When there are at least two symbols $\beta,\beta'$ in the bottom line, then as in the proof of Lemma~\ref{intersection:D}, there is a suspension data $\zeta''$ such that $|\zeta''_{\beta}|<\varepsilon$. And we can find paths in $D_{\pi,\varepsilon}$ that join $(\pi,\zeta)$ to $(\pi,\zeta'')$ and $(\pi,\zeta')$ to $(\pi,\zeta'')$. 
\end{itemize}

\end{proof}

\begin{proof}[Proof of Theorem~\ref{end}]
It is well known that $\HH(\emptyset)$ has only one topological end. The stratum $\QQ(-1,-1,-1,-1)$ can be naturally identified with $\HH(\emptyset)$, hence has also one topological end (see also \cite{B2}). 

For the other cases, we combine Theorem~\ref{extrc:cc},  Lemma~\ref{intersection:D} and~\ref{D:connected} to see that the set $\cup_{\pi\in C} Z(D_{\pi,\varepsilon})$, for $C$ the extended Rauzy class corresponding to $\CC$, is a subset of the same connected component of $\CC_{1,\varepsilon}$. 

By Lemma~\ref{rectangle}, any connected component of $\CC_{1,\varepsilon}$ intersects the subset $\cup_{\pi\in C} Z(D_{\pi,\varepsilon})$. Hence there is only one connected component in $\CC_{1,\varepsilon}$.
\end{proof}

\appendix
\section{Geometric description of a neigborhood of the boundary of a stratum}

In this this part, we present  natural splittings of a generic flat surface in $\mathcal{C}_{1,\varepsilon}$. We show some examples of the difficulties that can arise if one wants to relate two possible ``configurations'' that can occur in the same connected component of the corresponding strata of the moduli space of flat surfaces.

\subsection{Configurations of rigid collections of saddle connections}

Let $S$ be a surface in $\mathcal{C}_{1,\varepsilon}$. The set of lengths of saddle connections is discrete, hence, there exists a saddle connection whose length is minimal. 

Such saddle connection is not necessarily unique even if the surface is generic. Indeed, if two saddle connections on a translation surface are homologous, then they are necessary parallel and of the same length. This property is preserved by any small deformation of the surface inside the ambient stratum. 
One can show that the converse is true: if the ratio of the length of two saddle connections  is constant for any small deformation of a translation surface, then they are homologous. The analogous notion for half-translation surfaces is ``\^homo\-logous''. In this case, two \^homo\-logous saddle connections are parallel and the ratio of their lengths is in $\{1/2,1,2\}$ (see \cite{EMZ,Masur:Zorich} for more details).

Hence, a generic surface in  $\mathcal{C}_{1,\varepsilon}$ naturally defines a splitting of the surface by a collection of small and parallel saddle connections. This collection is preserved under small perturbation of the surface.
We will call \emph{configuration}, the discrete data associated to this splitting which is preserved under any small deformation of the surface (see also \cite{EMZ,Masur:Zorich}). Note that, contrary to the case of generalized permutations, there is no canonical way to describe all the configurations that can appear in a stratum, or in a connected component of a stratum (see \cite{Masur:Zorich,B1}).

\begin{Example}\label{ex:h211}
Consider an element $S_0\in \mathcal{H}(0)$ and $S_1\in \mathcal{H}(1,1)$. We slit each of these two surfaces along a segment of the same length and with the same direction such that exactly one endpoint of the segment in $S_1$ is a singularity.  We get two flat surfaces with one boundary component each, and each of these boundaries consists of two saddle connections. Now let $S$ be the translation surface obtained by gluing the two previous surfaces with boundary so that we get a closed translation surface. The segments of the boundary components of $S_1$ and $S_2$ correspond in $S$ to a pair $(\gamma_1,\gamma_2)$ of homologous saddle connections.  We see that the surface $S$ has three singularities: one corresponds to the singularity of $S_1$ that does not intersect the segment, and two correspond to the endpoints of the segments, hence, are of angle $2\pi+2\pi$ and $2\pi+4\pi$. Therefore $S$ is in the stratum $\mathcal{H}(2,1,1)$. The \emph{configuration} of $(\gamma_1,\gamma_2)$ can be seen as the combinatorial data associated to this construction.

Similarly, one can define a surface $S'$ with the same construction as before, but using a surface $S_2'$ in $\mathcal{H}(2)$ and with a segment that does not intersect the degree two singularity. It is easy to see that the surface $S'$ is also in the same stratum $\mathcal{H}(2,1,1)$. The corresponding saddle connections $(\gamma_1',\gamma_2')$ have a different configuration. 
\end{Example}

Note that a configuration of homologous saddle connections is also preserved by the $SL(2,\R)$ action on the stratum. Therefore, by ergodicity of this action, as soon as there exists a translation surface with a collection of homologous saddle connections, then such collection exists on almost all translation surfaces of the same connected component of stratum. Moreover, a result of Eskin and Masur asserts that the number of collections of saddle connections realizing a given configuration on a generic surface has quadratic asymptotics (see \cite{EM}). However, there is no lower bound on the length of a saddle connection that would realize a given configuration.

In the case of the surfaces $S,S'$ of Example~\ref{ex:h211}, by choosing sufficiently small segments in the construction, we can assume that the saddle connections $\gamma_i, \gamma_i'$ are very small and the surfaces are in the $\varepsilon$-boundary of the stratum. There exists a pair of homologous saddle connections $(\gamma_1'',\gamma_2'')$ on $S$ that realizes the configuration of $(\gamma_1',\gamma_2')$, but the saddle connections $\gamma_i''$ might be very long, and therefore, it can be difficult to shrink $\gamma_i''$ by staying in the $\varepsilon$-boundary of the stratum.
A possible solution is to look at some other saddle connections, not too long, that would correspond to an intermediary configuration. In our case, it is easy to find on $S$ and on $S'$ a simple saddle connection joining for instance a singularity of degree 1 and a singularity of degree 2. Here by simple we mean that no other saddle connection is homologous to it.
However, such approach  depends on the geometry of $S$ and $S'$, since there does not exists an analogous to the Rauzy induction for configurations, \emph{i.e.}
 a canonical operation that relates all the configurations of a connected component of a stratum. Furthermore, as we will see, the relations between configurations and the geometry of the surface are not simple. 

\subsection{Some strange examples}\label{ex:qi9}
In this section we present two examples that correspond to the  stratum $\emph{Q}(-1,9)$. This stratum is nonconnected and the only known proof of this result uses extended Rauzy classes (see \cite{Lanneau:cc,Boissy:Lanneau}).

\subsubsection{First example} \label{Ex1}
Let $S\in \QQ(-1,9)$ with the following decomposition: there exists two closed saddle connections $\gamma_1$ and $\gamma_2$ that start and end at the singularity of order 9 and that are the boundary of a cylinder, and such that no other saddle connections are \^homo\-logous to the $\gamma_i$.
\begin{figure}[htb]
\begin{center}
\psfrag{1}{$\gamma_1$}
\psfrag{2}{$\gamma_2$}
\psfrag{k}{$k\pi$}
\includegraphics[scale=0.65]{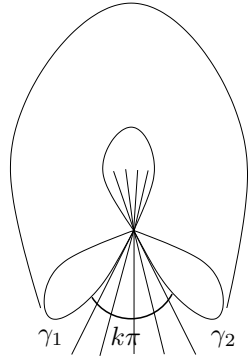}
\caption{Cylinder attached to a singularity, bounded by two saddle connections}
\end{center}
\end{figure}

The angle between $\gamma_1$ and $\gamma_2$ is $k\pi$, for $k\in \{1,2,3,4\}$. Lanneau (\cite{Lanneau:cc}) has proven the following.
\begin{itemize}
\item If $k\in \{1,2,4\}$ then $S$ belongs to the \emph{regular} connected component of $\QQ(-1,9)$, \emph{i.e.} there exists a simple saddle connection joining the pole to the zero of order 9.
\item If $k=3$, then $S$ belongs to the \emph{irreducible} connected component of $\QQ(-1,9)$, and hence, there is no simple saddle connection joining the pole to the zero of order 9.
\end{itemize}

Note that in the stratum $\QQ(-1,1+4n)$, for $n\geq 3$, one can find surfaces with analogous decomposition by a pair of saddle  connections bounding a cylinder. But in this case, for all parameters $k$, there exists a simple saddle connection joining the pole to the zero of order $1+4n$, because the stratum $\QQ(-1,1+4n)$ is connected.

\subsubsection{Second example}
Let $S\in \QQ(-1,9)$ and $\gamma$ a simple closed horizontal saddle connection of $S$ joining the zero to itself. It defines at the zero an unordered pair of angles $(\alpha \pi,\beta \pi)$ with $\alpha,\beta\in \mathbb{N}^* $ and $\alpha+\beta=11$. We assume that this pair is $(4\pi,7\pi)$. Then $S$ can be obtained from the stratum by the following construction: we start from a surface $S_0\in \QQ(-1,1,4)$ and  we denote by $P_1$ and $P_2$ the two singularities of order 1 and 4.  For each $P_i$, we choose an angular sector of angle $\pi$,  between two consecutive horizontal separatrices. We denote by $I$ and $II$ these sectors.
Then, we choose a path $\nu$ transverse to the horizontal foliation and  without self intersections that joins the sector $I$ of $P_1$ to the sector $II$ of $P_2$ (such path always exists by a result of Hubbard and Masur \cite{Hubbard:Masur}). Then, we cut the surface along this path and paste in a ``curvilinear annulus'' with two opposite sides isometric to $\nu$, and with horizontal sides of length $\varepsilon$. We get a surface with two boundary components and each of these components consists of a closed saddle connection.
These two saddle connections are parallel and of the same length, and therefore, we can isometrically glue them together and get a flat surface which is in $\QQ(-1,9)$, and with a saddle connection that realize the configuration described above. 

One can show that the connected component containing the surface obtained by this construction does not depend on the choice of the path $\nu$ (see \cite{B2}, Lemma ~4.5), once fixed the pair of sectors $I$ and $II$. Also, the resulting connected component does not change if we change continuously the initial surface. However, since the numbers of possible sectors (respectively $3$ and $6$) are not relatively prime, there remains a parameter in $\mathbb{Z}/3\mathbb{Z}$. 

In this precise case, we have the following result: fix the sector $II$. Among the 3 possible choices for sector $I$, the following possibilities hold:
\begin{itemize}
\item Two choices we give a surface in the \emph{regular} connected component of $\QQ(-1,9)$, and hence, on the resulting flat surface, there exist a simple saddle connection joining the pole to the zero of order~9.
\item The last choice leads to the \emph{irreducible} connected component of $\QQ(-1,9)$ and hence, on the resulting flat surface, we cannot find a simple saddle connection joining the pole to the zero of order~9.
\end{itemize}

\begin{figure}[htb]
\begin{center}
\psfrag{a}{a)}
\psfrag{b}{b)}
\psfrag{c}{c)}
\psfrag{1}{\tiny 1}
\psfrag{2}{\tiny 2}
\psfrag{0}{\tiny 0}
\psfrag{3}{\tiny 3}
\psfrag{4}{\tiny 4}
\psfrag{5}{\tiny 5}
\psfrag{6}{\tiny 6}
\psfrag{2a}{}
\psfrag{2b}{}
\psfrag{2c}{}
\includegraphics[scale=0.6]{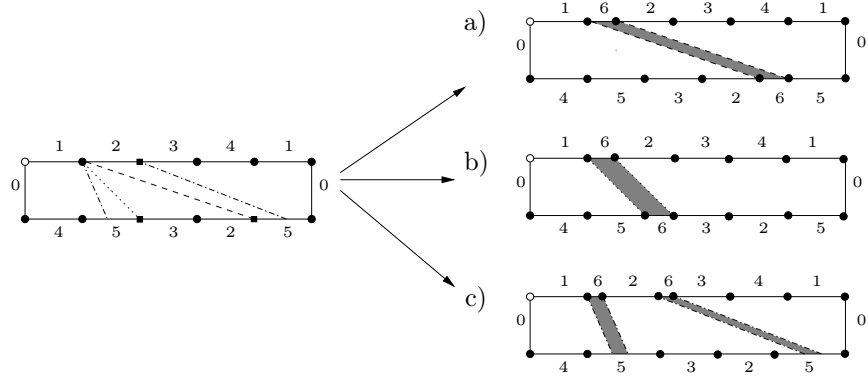}
\caption{Three similar surgeries of a surface in $\QQ(-1,4,1)$}
\label{Qi14rect}
\end{center}
\end{figure}

To prove this, we start from the flat surfaces obtained starting from a rectangle and gluing by translation or half turn the sides according to the left drawing of Figure~\ref{Qi14rect}. Then we perform the previous construction using the 3 possible paths that are indicated by dotted lines. This leads to 3 flat surfaces $a)$, $b)$ and $c)$.

\begin{figure}[htb]
\begin{flushleft}
\psfrag{a}{a)}
\psfrag{b}{b)}
\psfrag{c}{c)}
\psfrag{1}{\tiny 1}
\psfrag{2}{\tiny 2}
\psfrag{0}{\tiny 0}
\psfrag{3}{\tiny 3}
\psfrag{4}{\tiny 4}
\psfrag{5}{\tiny 5}
\psfrag{6}{\tiny 6}
\psfrag{irr}{ $\QQ_{irr}(-1,9)$}
\psfrag{reg}{$\QQ_{reg}(-1,9)$}
\includegraphics[scale=0.45]{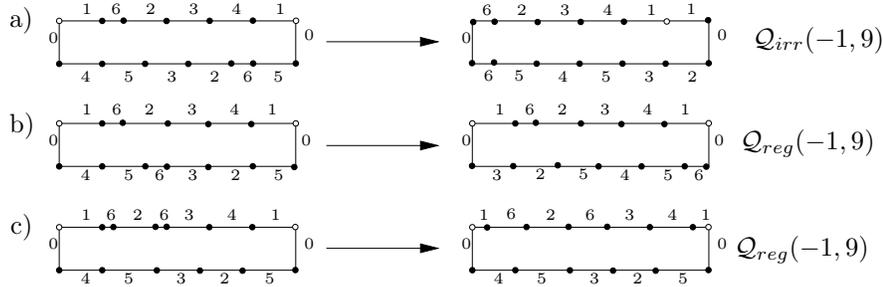}
\caption{Connected components of the surfaces of Figure~\ref{Qi14rect}}
\label{Qi14-cc}
\end{flushleft}
\end{figure}

Continuously deforming these 3 surfaces and after some cutting and gluing, we get three new surfaces in the same connected component as the initial ones. The new surfaces are shown on Figure~\ref{Qi14-cc}.

Then we have the following:
\begin{itemize}
\item In the case $a)$, there is a flat cylinder corresponding to the closed vertical geodesics starting from the segment~6. According to the criteria of Lanneau (see \cite{Lanneau:cc}, section 5), this cylinder is simple. We can check that the angle between the corresponding two saddle connections is $3\pi$. Hence, by the criteria of Lanneau described in Example~\ref{Ex1}, the flat surface is in the irreducible connected component of $\QQ(-1,9)$.
\item In the case $b)$, the criteria of Lanneau says that the vertical saddle connection corresponding to the segment labelled ``0'' is simple, hence we are in the regular connected component of $\QQ(-1,9)$.
\item  In the case $c)$, it is clear that the saddle connection corresponding to the segment labelled ``6'' is simple, hence we are in the regular connected component of $\QQ(-1,9)$.
\end{itemize}

\section*{Acknowledgments}
I  thank Anton Zorich, Pascal Hubert and Erwan Lanneau for encouraging me to write this paper, and for many discussions. I also thank the anonymous referee and Howard Masur for comments and remarks.

\end{document}